\documentclass[preprint,12pt]{elsarticle}

\usepackage{amssymb}
\usepackage{amsmath}
\usepackage{amsthm}
\DeclareMathOperator{\Hessian}{H}
\providecommand{\bs}[1]{{\mathbf{#1}}}

\newtheorem*{remark}{Remark}

\journal{}

\begin{document}
	
	\begin{frontmatter}
		
		
		
		\title{Multinode Shepard collocation method for pricing of financial derivatives}
		
		
		\author[address-CS]{Francesco Dell'Accio}
		\ead{francesco.dellaccio@unical.it}
		
		\author[address-CS]{Filomena Di Tommaso\corref{corrauthor}}
		\cortext[corrauthor]{Corresponding author}
		\ead{filomena.ditommaso@unical.it}

		\author[address-PA]{Elisa Francomano}
		\ead{elisa.francomano@unipa.it}
		
		\author[address-CS]{Clara Lorenzi}
		\ead{claralorenzi79@gmail.com}
		
		\address[address-CS]{Department of Mathematics and Computer Science, University of Calabria, Rende (CS), Italy}
		\address[address-PA]{Department of Engineering, University of Palermo,  Palermo, Italy}
		
		\begin{abstract}
	This paper explores the use of the multinode Shepard method for the numerical solution of the two-dimensional Black-Scholes equation. The proposed approach integrates a spatial approximation via the multinode Shepard operator with a temporal discretization based on the Backward Difference Formula. Numerical experiments are presented to demonstrate the accuracy and effectiveness of the method.
		\end{abstract}

		
	\end{frontmatter}
	
	\section{Introduction}
	The financial markets are becoming increasingly complex, with not only equities but also numerous types of derivative financial instruments being traded. The market requires up-to-date information on the values of these derivatives every second of every day. This leads to a huge demand for fast and accurate computer simulations \cite{MILOVANOVIC2020205}. 
	Among the best-known and most widely used models is the Black-Scholes equation, introduced in 1973 by Fischer Black and Myron Scholes, which revolutionised quantitative finance by providing a theoretical framework for option pricing \cite{bs}.
	Although the original formulation of the Black-Scholes equation is one-dimensional, the increasing complexity of financial products has necessitated the development of more complex models capable of describing multivariate phenomena \cite{IKAMARI2020e00564}.
	One of the most widely used approaches to option pricing is Finite Difference (FD) method \cite{book}. This method achieves high accuracy and has a good convergence rate although it comes at a significant computational cost, which becomes particularly problematic when dealing with multiple assets. Efficient standard numerical methods for the PDE formulation include adaptive FD \cite{persson2007pricing, LOTSTEDT20071159,linde2009highly}, high-order compact schemes \cite{DURING20124462,parabolicproblem}, Alternating Direction Implicit (ADI) schemes \cite{Hout2010303,haentjens2012alternating}, Radial Basis Function (RBF) approximation \cite{LARSSON2008175,PETTERSSON200882}, Radial Basis Function Partition of Unity (RBF-PU) method \cite{safdari2015radial, shcherbakov2016radial, SHCHERBAKOV2016185}, and Radial Basis Function generated Finite Differences (RBF-FD) method \cite{MILOVANOVIC20181462, pfdrbf}. In \cite{articlebench} and \cite{von2019benchop} several methods for pricing of options are implemented and evaluated. In the paper \cite{PETTERSSON200882} it was shown that the RBF method compares well with an adaptive FD method. In \cite{SHCHERBAKOV2016185}, a RBF-PU method is introduced to overcome the primary limitation of global RBF approaches - namely, the dense and ill-conditioned nature of the collocation linear system. Furthermore, since the partition of unity framework naturally lends itself to parallelization, it significantly improves the overall computational time complexity.
	
	The aim of this paper is to apply the multinode Shepard method - a meshless Partition of Unity technique based on weighted interpolation \cite{dell2019rate}, which has proven effective in solving numerical problems in multidimensional settings \cite{DellAccio_Solving,DellAccioMelecon,DELLACCIO2024115896, Dell'Accio202487} - to the solution of the Black-Scholes equation. 
	In the numerical treatment of the Black-Scholes equation, particular attention has been paid to the temporal discretisation, a crucial element in guaranteeing the stability and accuracy of the solution. To this end, in literature \cite{ MILOVANOVIC2020205,PETTERSSON200882}, the Backward Differentiation Formula (BDF) of second order has been used. This formula, which belongs to the class of implicit multistep methods, known for their stability and ability to handle stiff problems, allows a good approximation of the time derivative to be obtained even in the presence of steep gradients or complex solution behaviour \cite{stabilitybdf}.
	
	The paper is organized as follows.
	In Section \ref{sec2} we recall the multinode Shepard method from a theoretical point of view, focusing on its definitions and properties. 
	In Section \ref{sec3} we present the Black-Scholes equation, describing its context and general formulation for the two-dimensional case.
	In Section \ref{sec4} we discuss the spatial approximation and in Section \ref{sec5} the temporal discretisation based on second order BDF methods.
	Finally, in Section \ref{sec6} we conduct various numerical experiments to asses the effectiveness of the numerical proposed method.
	
	\section{The multinode Shepard method}\label{sec2}
	Let $X= \{ \mathbf{x}_{1},\dots ,\mathbf{x}_{n}\}$ be a set of $n$ pairwise distinct points in a domain $\Omega \subset \mathbb{R}^s$ with associated data $f_i$ sampled from an unknown function $f: \Omega \to \mathbb{R}$, that is $f_i=f(\mathbf{x}_i)$, $i=1, \dots n$. Let $p \in \mathbb{N}$ and let us denote by $\mathbb{P}_p(\mathbb{R}^s)$ the space of polynomials in $s$ variables of total degree $p$, of dimension $\tau =\binom{s+p}{s}$. 

Let us assume that a covering $\mathcal{T}=\{t_1, \dots, t_m\}$ of the set $X$ is given, where each subset $t_j$ has cardinality $\tau$ and is unisolvent for interpolation by polynomials in $\mathbb{P}_p(\mathbb{R}^s)$.  
We set $t_j=\{ \mathbf{x}_{j_k} \mid k=1, \dots, \tau \}$ and denote by $p_j[f] \in \mathbb{P}_p(\mathbb{R}^s)$ the unique polynomial such that $p_j[f](\mathbf{x}_{j_k})=f(\mathbf{x}_{j_k})$.

The multinode Shepard operator, associated to the covering $\mathcal{T}$, is then defined by
	\begin{equation}
		\begin{array}{ccc}
			\mathcal{M}_{\mu }\left[ f\right] \left( \mathbf{x}\right)
			=\sum\limits_{j=1}^{m}B_{\mu ,j}\left( \mathbf{x}\right) p_{j}\left[ f\right]
			\left( \mathbf{x}\right) , & \mu >0, & \mathbf{x\in }\mathbb{R}^s,%
		\end{array}
		\label{multinode_operator}
	\end{equation}%
where
\begin{equation}\label{multinodebasisfunction}
		\begin{array}{cc}
			B_{\mu ,j}\left( \mathbf{x}\right) =\dfrac{\prod\limits_{\ell=1}^{\tau}\left\Vert 
				\mathbf{x}-\mathbf{x}_{j_{\ell}}\right\Vert ^{-\mu }}{\sum\limits_{k=1}^{m}%
				\prod\limits_{\ell=1}^{\tau}\left\Vert \mathbf{x}-\mathbf{x}_{k_{\ell}}\right\Vert
				^{-\mu }}, & j=1,\dots ,m,%
		\end{array}
	\end{equation}%
 are the multinode Shepard functions, generalizing the Little triangular Shepard functions \cite{dell2015approximation}. These functions satisfy the following properties:
	\begin{enumerate}
		\item $\sum\limits_{j=1}^m B_{\mu ,j}\left(\mathbf{x}\right) =1,$
		\item $B_{\mu ,j}\left( \mathbf{x}_{i}\right) =0$ for all $\mathbf{x}_{i}\notin t_j$,
		\item $\nabla B_{\mu ,j}\left( \mathbf{x}_{i}\right) =\mathbf{0},\text{ }\mu >1$,
		\item $\Hessian B_{\mu ,j}\left( \mathbf{x}_{i}\right) =\mathbf{0},\text{ }\mu >2$, 
		\item $\sum\limits_{j\in \mathcal{J}_{i}}B_{\mu ,j}\left( \mathbf{x}_{i}\right) =1$,
		\item $\sum\limits_{j\in \mathcal{J}_{i}}\nabla B_{\mu ,j}\left( \mathbf{x}%
		_{i}\right) =\mathbf{0}\text{, }\mu >1$,
		\item $\sum\limits_{j\in \mathcal{J}_{i}}\Hessian B_{\mu ,j}\left( \mathbf{x}%
		_{i}\right) =\mathbf{0}\text{, }\mu >2$,
	\end{enumerate}
	where $\nabla B_{\mu ,j}\left( \mathbf{x}\right)$ and $\Hessian B_{\mu ,j}\left( \mathbf{x}\right)$ are the gradient and the Hessian matrix of $B_{\mu ,j}(\mathbf{x})$ and $\mathcal{J}_{i}=\{j\in \{1,\dots ,m\}:i\in \{j_{1},\dots ,j_{\tau}\}\}$ is the set of indices of the $\tau$-tuples $t_j$ which contain the point $\mathbf{x}_{i}$.

In applying the multinode Shepard operator, it is convenient to express the local polynomial interpolants in \eqref{multinode_operator} in barycentric form
	\begin{equation}
		p_{j}\left[ f\right] \left( \mathbf{x}\right) =\sum\limits_{i=1}^{\tau}\lambda
		_{j,i}\left( \mathbf{x}\right) f_{j_{i}},
		\label{Lagrange_poly}
	\end{equation}
    where $\lambda
		_{j,i}\left( \mathbf{x}\right)$ are the fundamental Lagrange polynomials related to the subsets $t_j$, satisfying the condition
  \begin{equation}
\lambda_{j,i}\left( \mathbf{x}_{j_{k}}\right) =\left\{ 
\begin{array}{cc}
1, & i=k, \\ 
0, & \text{otherwise.}%
\end{array}%
\right.  \label{kronecker_delta_pro}
\end{equation}%
        After a rearrangement, it is possible to rewrite the operator $\mathcal{M}_{\mu }\left[ f\right]$ in terms of the function values 
	$$\mathcal{M}_{\mu }\left[ f\right] \left( \mathbf{x}\right)=\sum\limits_{i=1}^{n}\sum\limits_{j\in \mathcal{J}_{i}}B_{\mu ,j}\left( \mathbf{x}\right) \lambda _{j,i}\left( \mathbf{x}\right) f_{i}=\sum\limits_{i=1}^{n} W_{\mu,i}\left(\bs{x}\right) f_{i}.$$
	We notice that the functions 
	\begin{equation}\label{Wmuj}		W_{\mu,i}\left(\bs{x}\right)=\sum\limits_{j\in \mathcal{J}_{i}}B_{\mu ,j}\left( \mathbf{x}\right) \lambda _{j,i}\left( \mathbf{x}\right)
	\end{equation}
	constitute a cardinal basis and we refer to it as the multinode cardinal Shepard basis (for more details see \cite{DellAccio_Solving,DellAccioMelecon,DELLACCIO2024115896, Dell'Accio202487}).
	\section{The Black-Scholes problem}\label{sec3}
	From the mathematical point of view, the multi-dimensional Black-Scholes equation takes the form \cite{MILOVANOVIC2020205}
	\begin{equation}\label{problemabs}
		\left \{
		\begin{array}{l}
			\dfrac{\partial}{\partial \hat{t}}P\left(\hat{t}, \mathbf{x}\right)=\mathcal{L}P\left(\hat{t}, \mathbf{x}\right), \ \  \hat{t}\in \mathbb{R}_+, \ \mathbf{x}\in \mathbb{R}_+^s,\\
			P\left(0, \mathbf{x}\right)= \Phi \left(\mathbf{x}\right), \ \ \mathbf{x}\in \mathbb{R}_+^s,
		\end{array}
		\right.  
	\end{equation}
	where
	\begin{equation}
		\mathcal{L}P=r\sum\limits_{i=1}^s x_i \dfrac{\partial P}{\partial x_i}+ \dfrac{1}{2} \sum \limits_{i,j=1}^s {\sigma_i}{\sigma_j}\rho_{ij} x_ix_j \dfrac{\partial^2 P}{\partial x_i \partial x_j}- rP,
        \label{LP}
	\end{equation}
	with $P\left(\hat{t},\mathbf{x}\right)$ being the value of the option at time $\hat{t}$ when the underlying assets have the values specified by $\mathbf{x}=(x_1,\dots,x_s)$. The coefficient 
	$r$ in equation \eqref{LP} denotes the scaled short interest rate, $\sigma_i$ is the volatility of the asset $i$, $\rho_{ij}$  is the correlation between assets $i$ and $j$, and the upper limit $s$ of the sums indicates the number of underlying assets, which corresponds to the number of spatial dimensions in the problem. A common example of a contract function for a European basket call option is the average option, mathematically expressed as
	\begin{equation*}
		\Phi\left(\mathbf{x}\right)=\max \left(\dfrac{1}{s}\sum \limits_{i=1}^s x_i - {K},0\right), 
	\end{equation*}
	where $K \in \mathbb{R}$ is the strike price.
	The special case of the two-dimensional Black-Scholes equation, which is the objective of investigation of the present paper, is useful in pricing options that depend on two underlying assets or in scenarios where there is greater interaction between variables \cite{hull2021}. For simplicity, we set $\left(x_1,x_2\right)=\left(x,y\right)$ and equation \eqref{LP} becomes	
	\begin{equation}\label{bs2}
		\mathcal{L}P =\ r\left(x\dfrac{\partial P}{\partial x}+ y \dfrac{\partial P}{\partial y}\right) + \dfrac{1}{2}\left(\sigma_1^2 x^2 \dfrac{\partial^2 P}{\partial x^2} + \sigma_2^2 y^2 \dfrac{\partial^2 P}{\partial y^2}\right)
		+ \rho \sigma_1 \sigma_2 x y \dfrac{\partial^2 P}{\partial x \partial y} - rP,
	\end{equation}	
	where we assume $\rho_{11}=\rho_{22}=1$ and $\rho_{12}=\rho_{21}=\rho$ \cite{MILOVANOVIC2020205}.
	\subsection{Boundary conditions}
	Janson and Tysk in \cite{JANSON_TYSK_2006}  demonstrate that the problem under consideration is well-posed even without explicit boundary conditions, both in space and time, provided that the growth at infinity is suitably restricted. Consequently, we need to restrict the problem to a finite domain and impose only near-field and far-field boundary conditions \cite{MILOVANOVIC2020205,PETTERSSON200882}.
	The near-field boundary is fixed at the point $x=y=0$, where the following condition is applied
	\begin{equation} \label{condizioneP}
		P\left(0,0,t\right)=0, \ \ 0\leq t \leq T. 
	\end{equation}
	At the far-field boundary we utilize the asymptotic solution
	\begin{equation}\label{conP}
		P\left(x,y,t\right)=\frac{x+y}{2}-Ke^{\bar{r}t}, \ 0\leq t \leq T, \ \Vert \bs{x}\Vert \rightarrow \infty .
	\end{equation}

	\section{Spatial approximation}\label{sec4}
    Under the above assumptions, let $\Omega\subset \mathbb{R}^2$ be the triangle of vertices  $(0,0)$, $(8,0)$ and $(0,8)$. We assume that a node set $X=\{\bs{x}_1,\dots, \bs{x}_{n}\}$ is given together with a covering $\mathcal{T}$ and that the points in $X$ are organized such that $\bs{x}_1,\dots, \bs{x}_{n_I}$ are interior points of $\Omega$, $\bs{x}_{n_I+1}$ is the origin and $\bs{x}_{n_I+2},\dots, \bs{x}_n$ are points on the boundary line $x+y=8$.
    
	We assume that the solution of the problem \eqref{problemabs} subject to the boundary conditions \eqref{condizioneP} and \eqref{conP} is approximated, at a fixed time step $t$, by a multinode Shepard operator $\mathcal{M}_\mu [P;t](\bs{x})$. Consequently, by imposing the differential conditions \eqref{bs2}, we have
	\begin{equation}
		\begin{split}
			\mathcal{L}P\approx & \mathcal{LM}_\mu[P;t](\bs{x})\\&=  \left(x\dfrac{\partial \mathcal{M}_\mu}{\partial x}(\bs{x})+ y \dfrac{\partial \mathcal{M}_\mu}{\partial y}(\bs{x})\right)+ \dfrac{1}{2}\left(\sigma_1^2x^2\dfrac{\partial^2 \mathcal{M}_\mu}{\partial x^2}(\bs{x})+ \sigma_2^2y^2\dfrac{\partial^2 \mathcal{M}_\mu}{\partial y^2}(\bs{x})\right)\\
			& \quad +\rho\sigma_1\sigma_2 x y \dfrac{\partial^2 \mathcal{M}_\mu}{\partial x \partial y}(\bs{x})-r\mathcal{M}_\mu(\bs{x}),
		\end{split}
	\end{equation}
	where, in order to simplify notation, we use $\mathcal{M}_\mu$ instead of $\mathcal{M}_\mu[P;t](x,y)$.
	By computing the first and second order derivatives of the multinode cardinal basis functions \eqref{Wmuj}, it follows that
	\begin{equation}\label{eqderivate}
		\begin{split}
\mathcal{LM}_\mu\left(\bs{x}\right)=&\sum\limits_{i=1}^n\left(\sum\limits_{j \in \mathcal{J}_i}r\left(x\dfrac{\partial W_{\mu,j}}{\partial x}\left(\mathbf{x}\right)+ y \dfrac{\partial W_{\mu,j}}{\partial y}\left(\mathbf{x}\right)\right)\right. \\&\left. +\dfrac{1}{2}\left(\sigma_1^2x^2\dfrac{\partial^2 W_{\mu,j}}{\partial x^2}\left(\mathbf{x}\right)+ \sigma_2^2y^2\dfrac{\partial^2 W_{\mu,j}}{\partial y^2}\left(\mathbf{x}\right)\right)\right.\\&
			+\left.\rho\sigma_1\sigma_2 x y \dfrac{\partial^2 W_{\mu,j}}{\partial x \partial y}\left(\mathbf{x}\right)-rW_{\mu,j}\left(\mathbf{x}\right)\right)P_i(t),
		\end{split}
	\end{equation}
	where $P_i(t)=\mathcal{M}_{\mu}[P;t](\bs{x}_i)$ is the approximation of the solution $P(t,\bs{x})$ of the problem \eqref{problemabs} at the point $\mathbf{x}_i$. 
   
	By evaluating \eqref{eqderivate} at the interior points $\bs{x}_k=(x_k,y_k)$ and by separating the terms related to the internal points from the ones related to the boundary points, we get the spatial matrices
    \[
    \begin{array}{ccc}
    \mathcal{A}=\left[
    \begin{array}{cccc}
        \mathcal{A}_{1,1} & \dots & \mathcal{A}_{1,n_I} \\
        \mathcal{A}_{2,1} & \dots & \mathcal{A}_{2,n_I} \\
        \vdots & & \vdots \\
        \mathcal{A}_{n_I,1} & \dots & \mathcal{A}_{n_I,n_I}
    \end{array}
    \right],
    &
    &  \mathcal{B}=\left[
    \begin{array}{cccc}
        \mathcal{B}_{1,n_I+1} & \dots & \mathcal{B}_{1,n} \\
        \mathcal{B}_{2,n_I+1} & \dots & \mathcal{B}_{2,n} \\
        \vdots & & \vdots \\
        \mathcal{B}_{n_I,n_I+1} & \dots & \mathcal{B}_{n_I,n}
    \end{array}
    \right]
     \end{array}
    \]
with   
	\begin{equation}\label{matricespaziale}
		\begin{array}{ll}
			\mathcal{A}_{ki}=&\sum\limits_{j \in \mathcal{J}_i}r\left(x_k\dfrac{\partial W _{\mu,j}}{\partial x}\left(\mathbf{x}_k\right)+ y_k \dfrac{\partial W_{\mu,j}}{\partial y}\left(\mathbf{x}_k\right)\right) \\  &+
			\dfrac{1}{2}\left(\sigma_1^2x_k^2\dfrac{\partial^2 W_{\mu,j}}{\partial x^2}\left(\mathbf{x}_k\right)+ \sigma_2^2y_k^2\dfrac{\partial^2 W_{\mu,j}}{\partial y^2}\left(\mathbf{x}_k\right)\right)\\
			&+\rho\sigma_1\sigma_2 x_k y_k \dfrac{\partial^2 W_{\mu,j}}{\partial x \partial y}\left(\mathbf{x}_k\right)-rW_{\mu,j}\left(\mathbf{x}_k\right), \\
            &k=1,\dots,n_I, \, i=1,\dots,n_I,
		\end{array}
	\end{equation}
	
	\begin{equation}\label{matricespaziale2}
		\begin{array}{ll}
			\mathcal{B}_{ki}=&\sum\limits_{j \in \mathcal{J}_i}r\left(x_k\dfrac{\partial W _{\mu,j}}{\partial x}\left(\mathbf{x}_k\right)+ y_k \dfrac{\partial W_{\mu,j}}{\partial y}\left(\mathbf{x}_k\right)\right) \\  &+
			\dfrac{1}{2}\left(\sigma_1^2x_k^2\dfrac{\partial^2 W_{\mu,j}}{\partial x^2}\left(\mathbf{x}_k\right)+ \sigma_2^2y_k^2\dfrac{\partial^2 W_{\mu,j}}{\partial y^2}\left(\mathbf{x}_k\right)\right)\\
			&+\rho\sigma_1\sigma_2 x_k y_k \dfrac{\partial^2 W_{\mu,j}}{\partial x \partial y}\left(\mathbf{x}_k\right)-rW_{\mu,j}\left(\mathbf{x}_k\right), \\
            & k=1,\dots, n_I, \, i=n_I+1,\dots,n.
		\end{array}
	\end{equation}
	\section{Temporal discretization}\label{sec5}
	For the time discretization of the equation \eqref{bs2},  a Backward Differentiation Formula of order $2$ (BDF2) is adopted \cite{MILOVANOVIC2020205}. These time-stepping schemes belong to the class of implicit methods and are characterized by using the solution values computed in the previous $2$ time steps.
	
	To implement the BDF2 method, it is necessary to determine the initial solution values for the first time step using the Backward Euler method. This approach ensures a consistent scheme while preserving the stability and accuracy properties of the method starting from the given initial data \cite{stabilitybdf}. 
	To this goal, in line with \cite{MILOVANOVIC2020205}, we discretize the time interval $[0,T]$ into $M$ steps $\Delta t^\ell=t^\ell-t^{\ell-1}$, $\ell=1, \dots, M$ with $t_0=0$ and $t_M=T$. Moreover, we set 
	$$\bs{P}^\ell=[P^{\ell}_1,\dots P^{\ell}_{n_I},P^{\ell}_{n_I+1},\dots, P^{\ell}_n]^T,$$ with
	$$P^{\ell}_i=P_i(t^\ell)= \mathcal{M}_{\mu}[P;t^\ell](\bs{x}_i)\approx P(t^\ell,\bs{x}_i),\quad i=1,\dots,n,\, \ell=1,\dots,M.
	$$
The initial condition in equation \eqref{problemabs}, in its discrete form, is given by 
	\begin{equation}\label{u_0}
		\bs{P}^0=\left[\Phi\left({x_1}\right), \dots,\Phi\left({x_{n_I}}\right)\right]^T
	\end{equation}
and the time-stepping scheme for the initial time step is
	\begin{equation}
		\bs{P}^1-\bs{P}^0=\Delta t^1 \sum\limits_{i=1}^n\mathcal{L}W_{\mu,i}\left(\mathbf{x}\right)P^1_i,
	\end{equation}
that can be rewritten as
	\begin{equation}
		\bs{P}^1-\Delta t^1 \sum\limits_{i=1}^n\mathcal{L}W_{\mu,i}\left(\mathbf{x}\right)P^1_i=\bs{P}^0.
	\end{equation}
Partitioning the previous sum by distinguishing the $n_I$ internal nodes from the $n-n_I$ boundary nodes, the previous equation can be expressed as 
	\begin{equation}
		\bs{P}^1-\Delta t^1 \sum\limits_{i=1}^{n_I}\mathcal{L}W_{\mu,i}\left(\mathbf{x}\right){P}^1_i-\Delta t^1 \sum\limits_{i=n_I+2}^{n}\mathcal{L}W_{\mu,i}\left(\mathbf{x}\right){P}^1_i=\bs{P}^0, \label{timestep1} 
	\end{equation}
	where the lack of the term $\Delta t^1\mathcal{L}W_{\mu,n_I+1}\left(\mathbf{x}\right){P}^1_{n_I+1} $ is due to the presence of the null near-field boundary condition, as specified in equation \eqref{condizioneP}.
	By imposing the far-field boundary conditions \eqref{conP} directly into the equation \eqref{timestep1}, the contribution of the boundary points is explicitly incorporated in the second sum. This leads to a reformulation, where the values at the boundary are prescribed, ensuring consistency with the given conditions. As a result, the equation takes the following form
	\begin{equation}
		\bs{P}^1-\Delta t^1 \sum\limits_{i=1}^{n_I}\mathcal{L}W_{\mu,i}\left(\mathbf{x}\right){P}_i^1-\Delta t^1 \sum\limits_{i=n_I+2}^{n}\mathcal{L}W_{\mu,j}\left(\mathbf{x}\right)\left(\dfrac{x_i+y_i}{2}-Ke^{-rt^1}\right)=\bs{P}^0.
	\end{equation}
	Since the term $$\Delta t^1 \sum\limits_{i=n_I+2}^{n}\mathcal{L}W_{\mu,i}\left(\mathbf{x}\right)\left(\dfrac{x_i+y_i}{2}-Ke^{-rt^1}\right)$$
	is known, the previous equation becomes
	\begin{equation}
		\bs{P}^1-\Delta t^1 \sum\limits_{i=1}^{n_I}\mathcal{L}W_{\mu,j}\left(\mathbf{x}\right){P}_i^1=\bs{P}_i^0 +\Delta t^1 \sum\limits_{i=n_I+2}^{n}\mathcal{L}W_{\mu,i}\left(\mathbf{x}\right)\left(\dfrac{x_i+y_i}{2}-Ke^{-rt^1}\right).
	\end{equation}
	Evaluating at the collocation nodes $\bs{x}_k$, $k=1,\dots,n_I$, we obtain
	\begin{equation}\label{step1}
		\left(\mathbf{I}-\Delta t^1 \mathcal{A}\right)\bs{P}^1=\bs{P}^0+\Delta t^1 \mathcal{B}\, \mathbf{b},
	\end{equation}
	where $\mathbf{I}$ is the identity matrix of dimension $n_I \times n_I$, $\mathcal{A}$ is the spatial matrix for the interior nodes \eqref{matricespaziale}, $\mathcal{B}$ is the spatial matrix for the boundary nodes \eqref{matricespaziale2} and 
	\begin{equation}\label{termine_noto}
		\mathbf{b}=\left[\dfrac{x_i+y_i}{2}-Ke^{-rt^1}\right]_{i=n_I+2}^{n}.
	\end{equation}
	The time-stepping scheme progresses to the second time step, following the first step in the Backward Differentiation Formula (BDF). Specifically, for the second time step, the scheme is formulated as
	\begin{equation}
		\bs{P}^2-\dfrac{4}{3}\bs{P}^1+\dfrac{1}{3}\bs{P}^0=\dfrac{2}{3}\Delta t^2\sum\limits_{i=1}^n\mathcal{L}W_{\mu,i}\left(\mathbf{x}\right){P}_i^2.
	\end{equation}
	In line with the previous time step, by imposing the boundary conditions and by splitting the summation over the interior and boundary points and by  evaluating at the collocation nodes $\bs{x}_k$, we get 
	\begin{equation}\label{bdf2}
		\left(\mathbf{I}-\dfrac{2}{3}\Delta t^2 \mathcal{A}\right)\bs{P}^2=\dfrac{4}{3}\bs{P}^1-\dfrac{1}{3}\bs{P}^0+\dfrac{2}{3}\Delta t^2 \mathcal{B} \mathbf{b}.
	\end{equation}
	where $\mathbf{I}$, $\mathcal{A}$, $\mathcal{B}$ and $\mathbf{b}$ maintain the same definitions as in the previous time step. For the BDF2 scheme, the generic time step is then
	\begin{equation}
		\left(\mathbf{I}-\dfrac{2}{3}\Delta t^n \mathcal{A}\right)\bs{P}^n=\dfrac{4}{3}\bs{P}^{n-1}-\dfrac{1}{3}\bs{P}^{n-2}+\dfrac{2}{3}\Delta t^n \mathcal{B} \mathbf{b}.
	\end{equation}

	\section{Numerical results}\label{sec6}
	In this Section, we analyze various node configurations to identify node distributions capable of ensuring both high accuracy and numerical stability of the approximate solution computed through the multinode Shepard method. The experiments aim to evaluate both the mean absolute error and the maximum absolute error, thereby providing a solid foundation for the subsequent analysis.
	For each time step $t^\ell$, $\ell=1,\dots, M$, determine the mean absolute error using the formula 
	\begin{equation*}
		E_{mean}(t^{\ell})=\dfrac{1}{n_{eval}}\sum\limits_{i=1}^{n_{eval}} \vert \mathcal{M}_{4}[P;t^\ell](\boldsymbol{\xi}_i)-u^*(\boldsymbol{\xi}_i,t^\ell) \vert,
	\end{equation*}
	and the maximum absolute error using the formula 
	\begin{equation*} E_{\text{max}}(t^{\ell})= \max\limits_{1 \leq i \leq n_{eval}} \left|\mathcal{M}_{4}[P;t^\ell](\boldsymbol{\xi}_i) - u^*(\boldsymbol{\xi}_i,t^\ell) \right|,
	\end{equation*} 	
	where $u^*$ is a reference solution computed with a second order finite difference method on a very fine grid, $\mathcal{M}_{4}[P;t^\ell]$ is the approximate solution realized by the multinode Shepard operator and $\boldsymbol{\xi}_i,$ $ i=1,\dots, n_{eval}$ are the evaluation points arranged on a regular grid on the triangular domain.
    The multinode Shepard approximation is computed by using Halton nodes set \cite{halton}, Bézier node set \cite{cheney2009course}, Bézier and Waldron node set \cite{waldron}, Bézier and Waldron nodes set enriched with nodes lying on special straight lines. In all the experiments, the degree of the local polynomial interpolant is set to $p=2$ and the number of time step $M$ is fixed to $20$. The covering $\mathcal{T}$ is realized by selecting, for each collocation point $\bs{x}_i$, $i=1,\dots,n$, subsets of $6$ Leja points among $6+q$ nearest points by using the algorithm proposed in \cite{bos2010computing}. In all the experiments, we set $q=10$, except for those involving the Bézier and Waldron nodes, described in Section \ref{sec:bezier+waldron}, where we set $q=30$.

\begin{figure}
		\centering
		\parbox{0.32\linewidth}{\centering
			\includegraphics[width=1\linewidth]{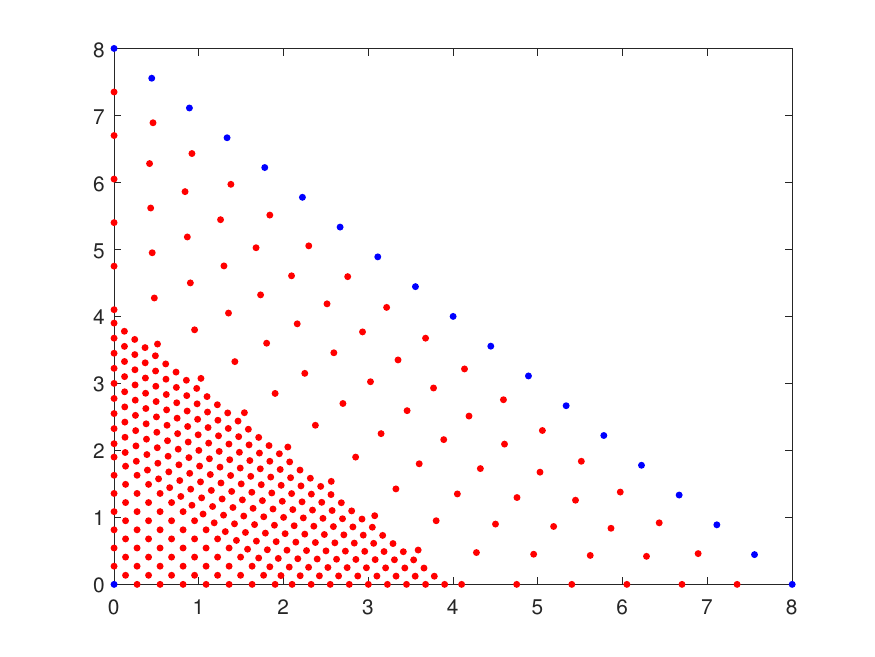}\\[-0.5ex]
            (a)
		}
        \parbox{0.32\linewidth}{\centering
			\includegraphics[width=1\linewidth]{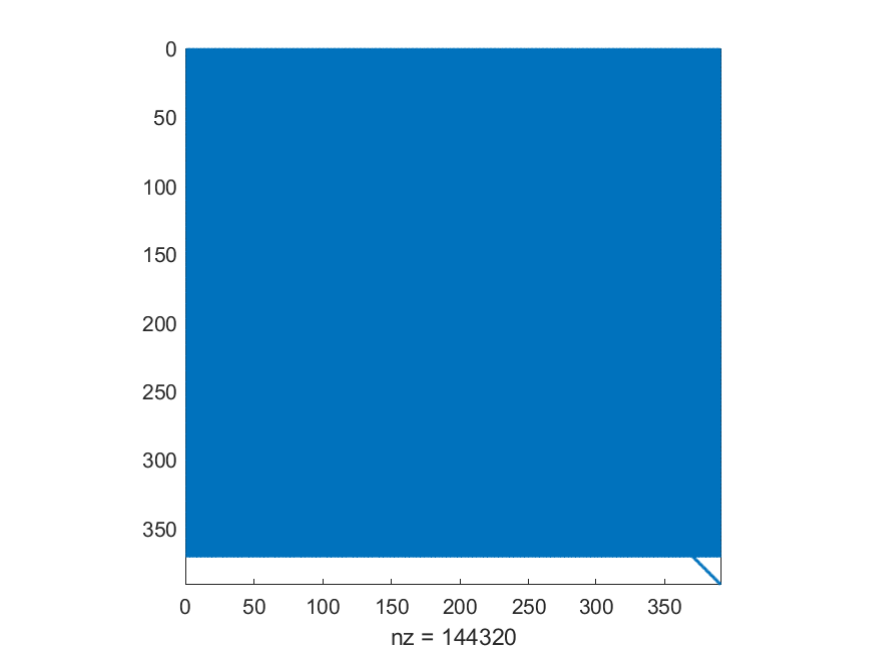}\\[-0.5ex]
            (b)
		}
		\hfill
		\parbox{0.32\linewidth}{\centering
			\includegraphics[width=1\linewidth]{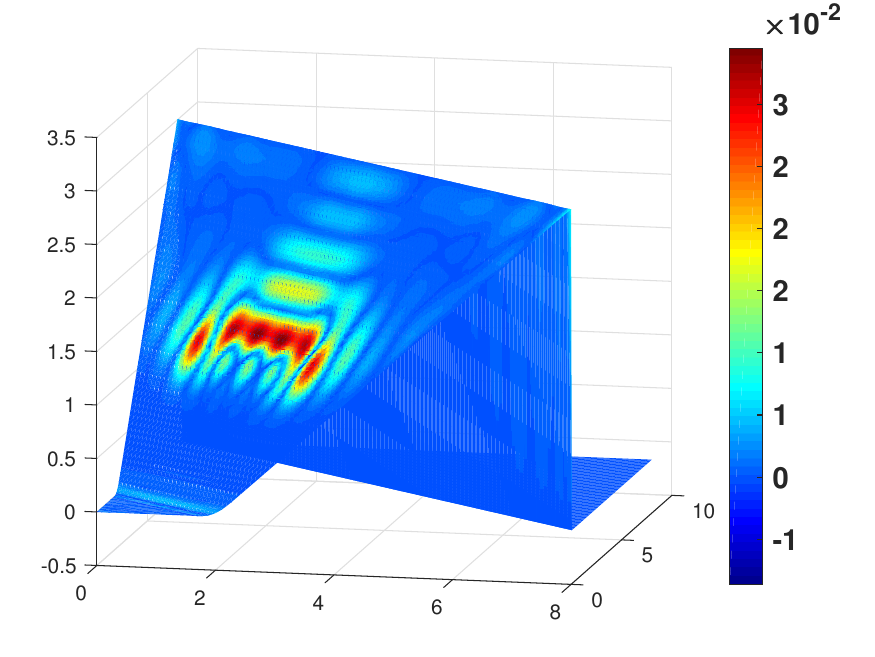}\\[-0.5ex]
            (c)
		}
        \caption{(a) Set of $370$ points coupled with $20$ boundary points. (b) Sparsity pattern of the collocation matrix related to node distribution. (c) Error on the surface of the approximated solution at the last time step.}
        \label{figlars}
	\end{figure}

    We compare the approximation accuracy of the Multinode Shepard method combined with the BDF2 temporal discretization scheme (MS-FD) with that obtained by the method proposed in \cite{LARSSON2008175}, based on multiquadric Radial Basis Function combined with the same temporal discretization (RBF-FD). For the latter, we make use of the non-uniform distribution of nodes in Figure \ref{figlars} (left), for which we show, in the same Figure, the sparsity pattern of the collocation matrix (center) and the error on the surface of the approximate solution (right). As stated by the authors in the referring paper \cite{LARSSON2008175}, such a distribution of nodes is chosen because it gives a clear improvement compared with a uniform distribution. 
 
	\subsection{Halton node set}
\begin{figure}
		\centering
		\parbox{0.32\linewidth}{\centering
		\includegraphics[width=\linewidth]{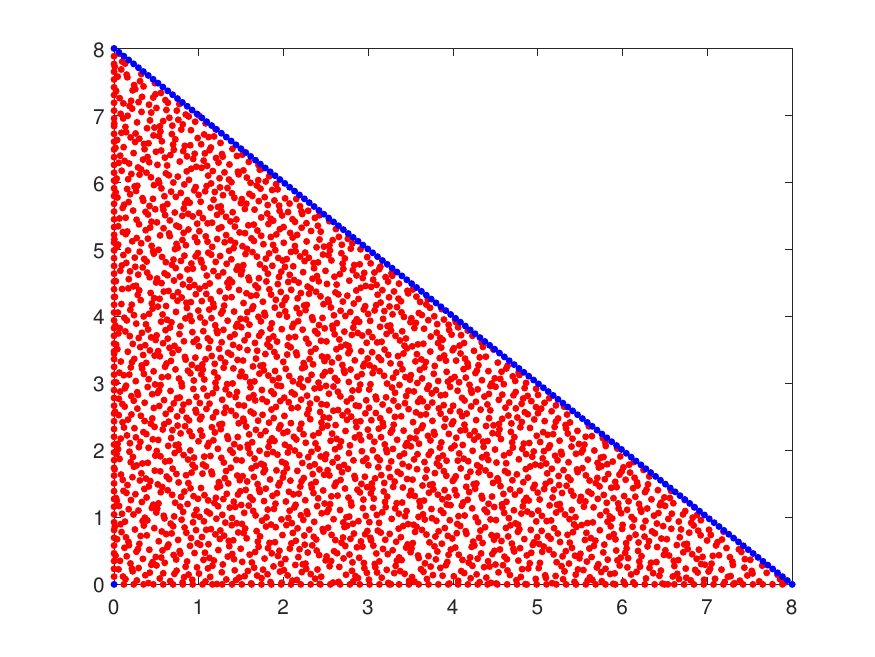}\\[-0.5ex]
		(a)
		}
		\parbox{0.32\linewidth}{\centering
			\includegraphics[width=\linewidth]{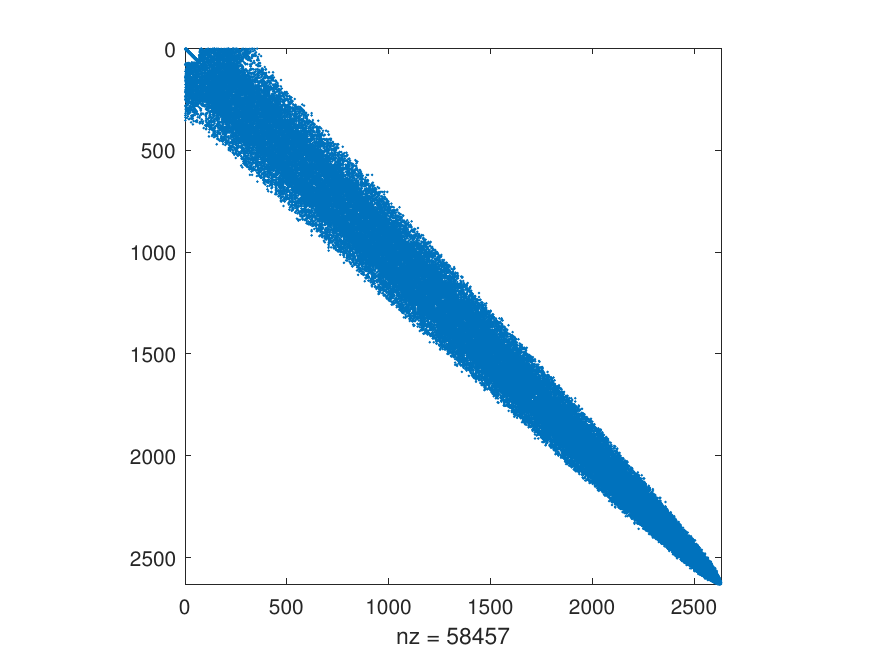}\\[-0.5ex]
			(b)
		}
		\parbox{0.32\linewidth}{\centering
			\includegraphics[width=\linewidth]{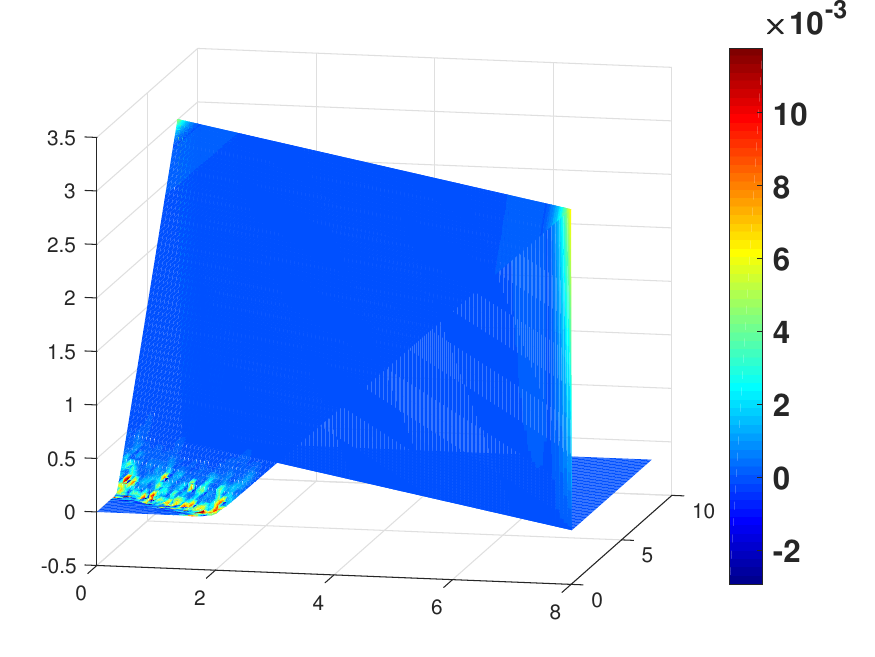}\\[-0.5ex]
			(c)
		}
		\caption{(a) Set of $2631$ Halton points coupled with $141$ boundary points. (b) Sparsity pattern of the collocation matrix related to the Halton points distribution. (c) Error on the surface of the approximated solution at the last time step with the MS-FD method.}
		\label{fig:haltonconfiguration}
	\end{figure}

	For the numerical analysis that follows, we have chosen a distribution of nodes that guarantees good coverage of the domain. For this purpose, we use quasi-random points, known for their uniform space-filling properties. In particular, we choose the well-known Halton sequence \cite{halton}.
	More precisely, we consider a set of $5000$ Halton points in $(0,1)\times (0,1)$ mapped linearly in the square $(0,8)\times (0,8)$ from which we extract those lying in the triangular domain defined by the vertices $(0,0)$, $(0,8)$ and $(8,0)$ and we couple this set with the origin of the axis and with $141$ equispaced nodes on the line $x+y=8$. In Figure \ref{fig:haltonconfiguration}, we show the node distribution, the sparsity pattern of the collocation matrix, and the pointwise error on the surface of the approximate solution realized by the MS-FD method. The condition number of the collocation matrix of the MS-FD method is about $8.4$ at the first time step, and $5.5$ from the second time step onward.
	In Table \ref{errore:halton}, we report the mean and max absolute errors on the various time steps by using MS-FD approximations as well as RBF-FD approximations.
    
    We observe, from Table \ref{errore:halton},  that the mean absolute errors remain stable as the number of time steps increases, suggesting good temporal stability in the proposed methods. The maximum error for MS-FD remains comparable to the RBF-FD maximum error. MS-FD outperforms the RBF-FD based approach in terms of mean absolute error by approximately one order of magnitude. Overall, the same table illustrates that MS-FD method provides superior accuracy compared to the RBF-FD scheme while maintaining stability and reliability over multiple time steps.

	\begin{table}
		\centering
		\footnotesize
		\begin{tabular}{c c c c c }
			\hline
			Time Step &$E_{mean}^{RBF-FD}$& $E_{max}^{RBF-FD}$  &$E_{mean}^{MS-FD}$&$E_{max}^{MS-FD}$ \\
			\hline 
			1&2.1679e-3&4.4602e-2&  1.9585e-4&  1.4363e-2\\
			2&2.1658e-3&4.5045e-2&  1.9507e-4&  1.3470e-2\\
			3&2.1554e-3&4.5312e-2&  1.9405e-4&  1.2691e-2\\
			4&2.1393e-3&4.5404e-2&  1.9333e-4&  1.2006e-2\\
			5&2.1182e-3&4.5335e-2&  1.9298e-4&  1.1401e-2\\
			6&2.0927e-3&4.5120e-2&  1.9270e-4&  1.0863e-2\\
			7&2.0638e-3&4.4778e-2&  1.9262e-4&  1.0384e-2\\
			8&2.0314e-3&4.4322e-2&  1.9245e-4&  1.0054e-2\\
			9&1.9965e-3&4.3767e-2&  1.9229e-4&  1.0057e-2\\
			10&1.9598e-3&4.3127e-2&  1.9193e-4&  1.0212e-2\\
			11&1.9220e-3&4.2414e-2&  1.9166e-4&  1.0333e-2\\
			12&1.8839e-3&4.1639e-2&  1.9134e-4&  1.0424e-2\\
			13&1.8455e-3&4.0811e-2&  1.9130e-4&  1.0531e-2\\
			14&1.8071e-3&3.9939e-2&  1.9151e-4&  1.0690e-2\\
			15&1.7695e-3&3.9033e-2&  1.9181e-4&  1.0822e-2\\
			16&1.7322e-3&3.8098e-2&  1.9244e-4&  1.0929e-2\\
			17&1.6956e-3&3.7143e-2&  1.9321e-4&  1.1014e-2\\
			18&1.6598e-3&3.6172e-2&  1.9419e-4&  1.1163e-2\\
			19&1.6246e-3&3.5191e-2&  1.9541e-4&  1.1345e-2\\
			20&1.5908e-3&3.4204e-2&  1.9690e-4&  1.1510e-2\\
			\hline
		\end{tabular}
		\caption{Mean and maximum absolute error for Halton node set shown in Figure \ref{fig:haltonconfiguration} on the various time steps by using the RBF-FD and MS-FD approximations.}
		\label{errore:halton}
	\end{table}
	
\subsection{Uniform distribution of nodes}
\begin{figure}
		\centering		
		\parbox{0.32\linewidth}{\centering
			\includegraphics[width=\linewidth]{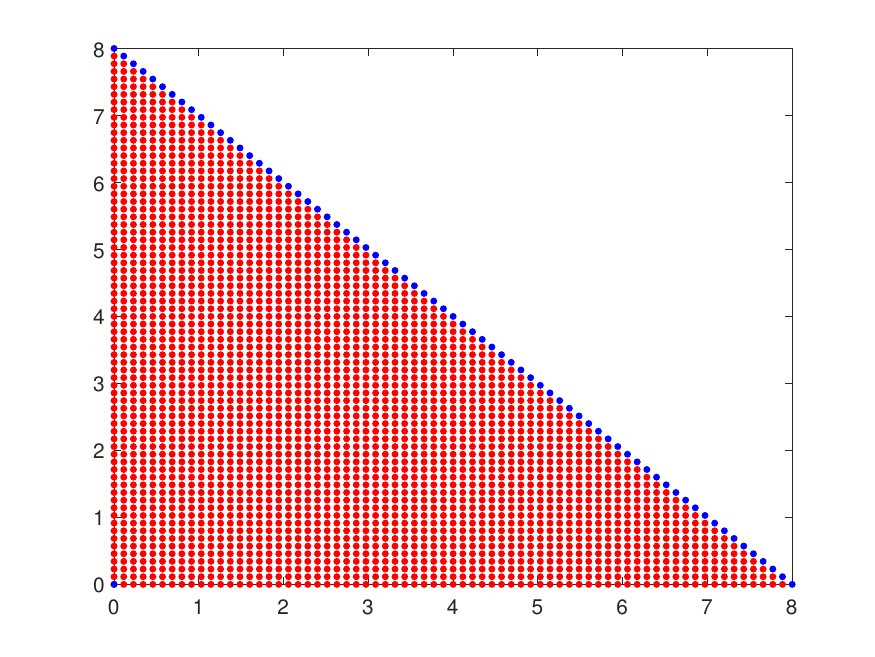}\\[-0.5ex]
			(a)
		}
		\parbox{0.32\linewidth}{\centering
			\includegraphics[width=1\linewidth]{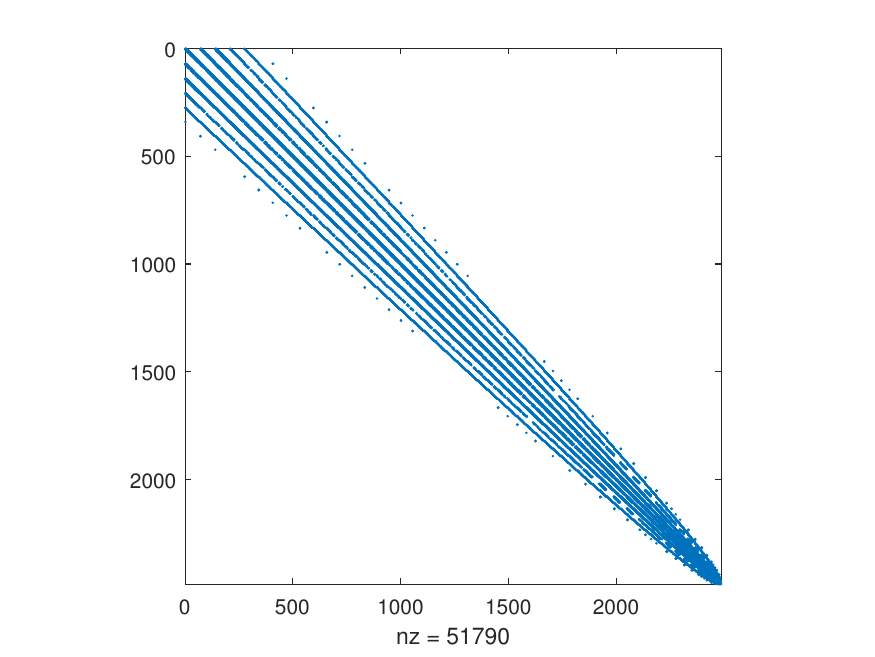}\\[-0.5ex]
			(b)
		}
		\parbox{0.32\linewidth}{\centering
			\includegraphics[width=1\linewidth]{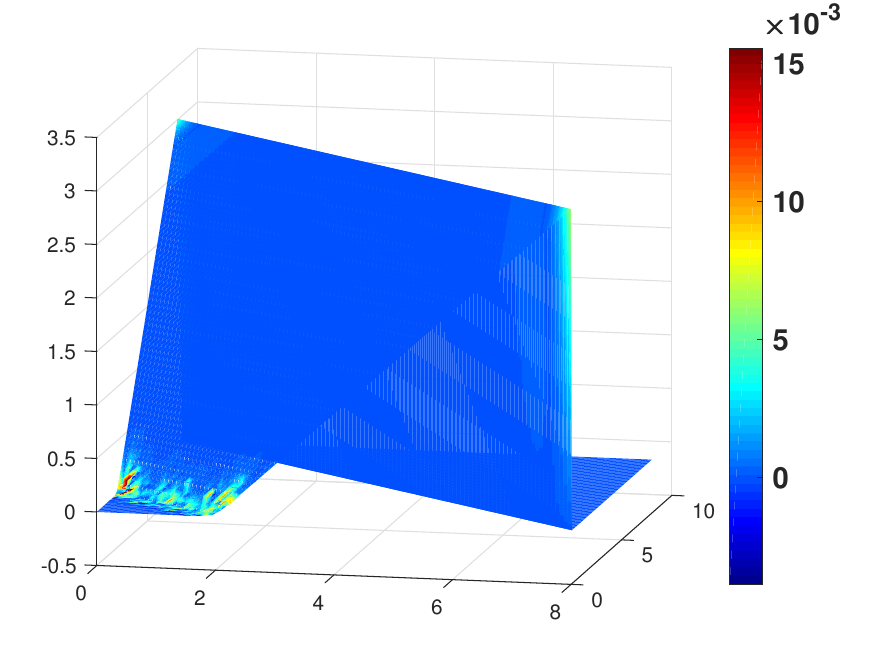}\\[-0.5ex]
			(c)
		}
		
		\caption{(a) Set of $2845$ uniformly distributed nodes coupled with $72$ boundary points (b) Sparsity pattern of the collocation matrix. (c) Error on the surface of the approximated solution at the last time step with the MS-FD method. }
		\label{fig:bezierconfiguration}
	\end{figure}

	To further investigate the performance of the MS-FD method under different spatial distribution of points, we now turn our attention to the uniform distribution of $\binom{n
    +2}{2}$ points $\left\{\left(\frac{i}{n},\frac{j}{n}\right): i,j \in \mathbb{N},\, i+j\leq n\right\}$  on the standard simplex.  This configuration of nodes is unisolvent for polynomial interpolation of total degree $n$ \cite[Theorem 4]{cheney2009course}. We refer to configurations of nodes that share this property as $n$-degree distributions. 
	In the experiments, we consider the $70$-degree uniform distribution on the triangle $(0,0)$, $(8,0)$ and $(0,8)$, which is obtained as the image of the $70$-degree uniform distribution on the standard simplex through the affine map $(x,y)\rightarrow (8x,8y)$. In Figure \ref{fig:bezierconfiguration}, we show the node distribution, the sparsity pattern of the collocation matrix, and the error on the surface of the approximated solution realized by the MS-FD method. 
	The condition number of the collocation matrix of the MS-FD method is about $4.1$ at the first time step, and $2.8$ from the second time step onward. In Table \ref{errore:bezier} we report the mean and maximum absolute error for the various time steps by using MS-FD approximations.

	Compared to previous results with Halton nodes in Table \ref{errore:halton} the MS-FD method exhibits a small increase in maximum absolute errors at each time step.
	The mean absolute errors, however, remain of the same order of magnitude to the Halton-based results, showing that the $70$-degree uniform distribution of nodes still produces stable and accurate solutions when used with the MS-FD method.
    
	\begin{table}
		\centering
		\footnotesize
		\begin{tabular}{c c c c c }
			\hline
			Time Step &$E_{mean}^{RBF-FD}$& $E_{max}^{RBF-FD}$  &$E_{mean}^{MS-FD}$&$E_{max}^{MS-FD}$ \\
			\hline
			1&2.1679e-3&4.4602e-2&  2.0861e-4&  1.9144e-2\\
			2&2.1658e-3&4.5045e-2&  2.0963e-4&  1.8376e-2\\
			3&2.1554e-3&4.5312e-2&  2.0954e-4&  1.7730e-2\\
			4&2.1393e-3&4.5404e-2&  2.0954e-4&  1.7183e-2\\
			5&2.1182e-3&4.5335e-2&  2.0972e-4&  1.6719e-2\\
			6&2.0927e-3&4.5120e-2&  2.0980e-4&  1.6322e-2\\
			7&2.0638e-3&4.4778e-2&  2.0976e-4&  1.5981e-2\\
			8&2.0314e-3&4.4322e-2&  2.0963e-4&  1.5687e-2\\
			9&1.9965e-3&4.3767e-2&  2.0943e-4&  1.5432e-2\\
			10& 1.9598e-3&4.3127e-2&  2.0905e-4&  1.5209e-2\\
			11&1.9220e-3&4.2414e-2&  2.0874e-4&  1.5013e-2\\
			12&1.8839e-3&4.1639e-2&  2.0842e-4&  1.4838e-2\\
			13&1.8455e-3&4.0811e-2&  2.0845e-4&  1.4683e-2\\
			14&1.8071e-3&3.9939e-2&  2.0883e-4&  1.4543e-2\\
			15&1.7695e-3&3.9033e-2&  2.0940e-4&  1.4415e-2\\
			16&1.7322e-3&3.8098e-2&  2.1025e-4&  1.4298e-2\\
			17&1.6956e-3&3.7143e-2&  2.1127e-4&  1.4190e-2\\
			18&1.6598e-3&3.6172e-2&  2.1250e-4&  1.4526e-2\\
			19&1.6246e-3&3.5191e-2&  2.1403e-4&  1.4861e-2\\
			20&1.5908e-3&3.4204e-2&  2.1595e-4&  1.5185e-2\\
			
			\hline
		\end{tabular}
		\caption{Mean and maximum absolute error for the uniform distribution of nodes in Figure \ref{fig:bezierconfiguration} on the various time steps by using the MS-FD method.}
		\label{errore:bezier}
	\end{table}

	\subsection{Uniform distribution of nodes combined with Waldron nodes}\label{sec:bezier+waldron}
\begin{figure}
		\centering
		\parbox{0.32\linewidth}{\centering
			\includegraphics[width=\linewidth]{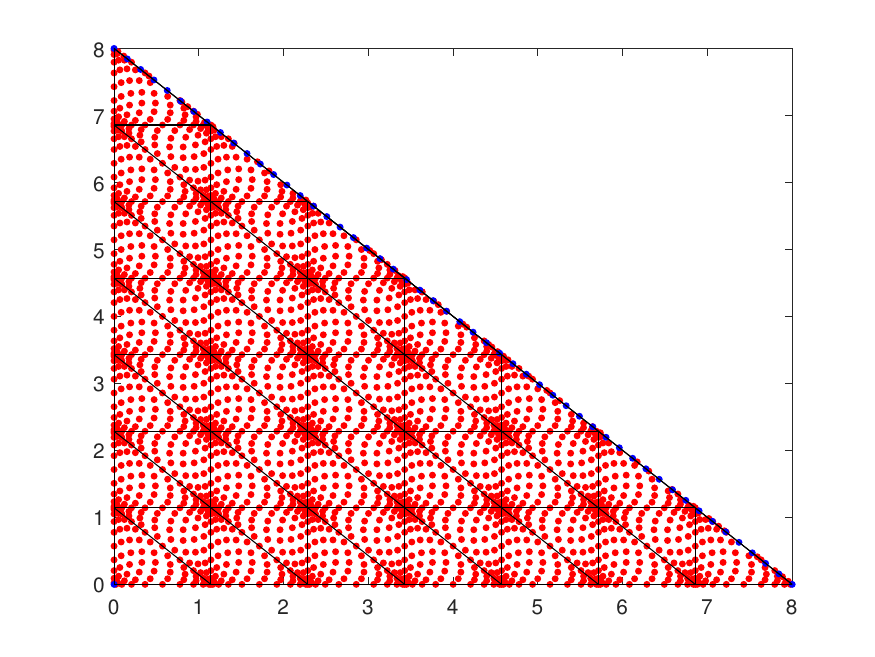}\\[-0.5ex]
			(a)
		}
		\parbox{0.32\linewidth}{\centering
			\includegraphics[width=1\linewidth]{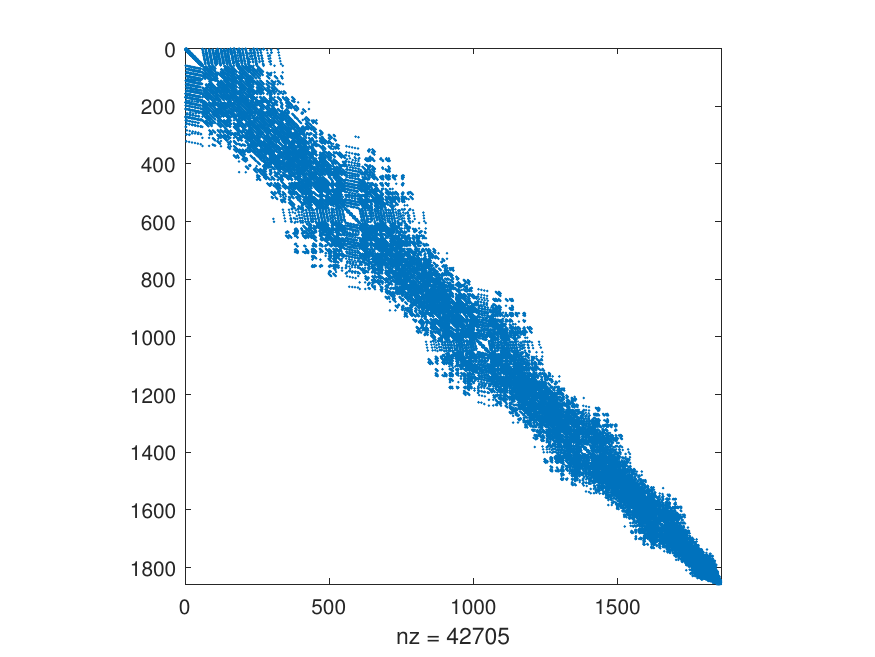}\\[-0.5ex]
			(b)
		}
		\parbox{0.32\linewidth}{\centering
			\includegraphics[width=1\linewidth]{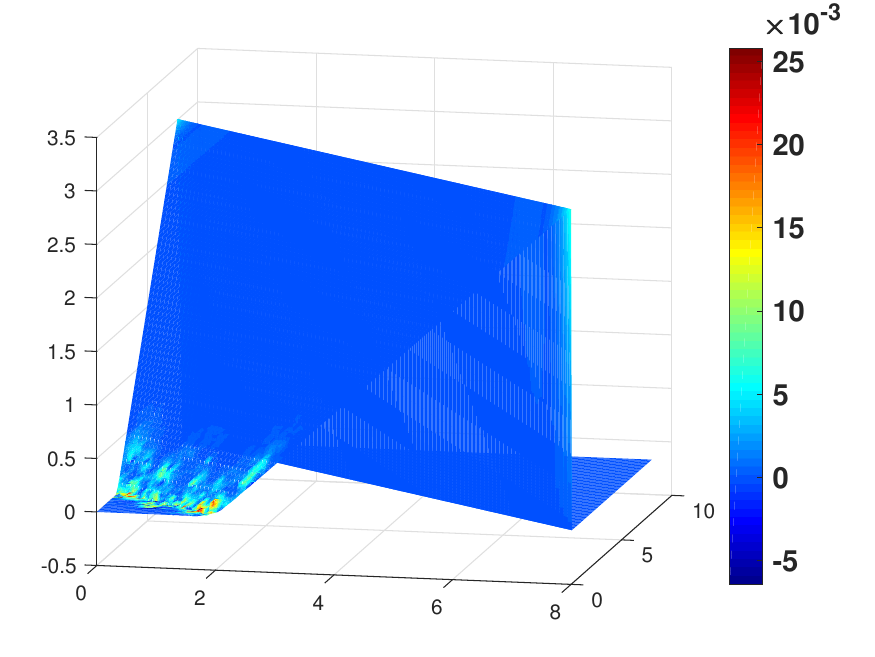}\\[-0.5ex]
			(c)
		}
		
		\caption{(a) Set of $1774$ uniform distribution of nodes combined with Waldron nodes coupled with $146$ boundary points. (b) Sparsity pattern of the collocation matrix. (c) Error on the surface of the approximated solution at the last time step with the MS-FD method.}
		\label{fig:bezierwaldronconfiguration}
	\end{figure}

	In this section, we explore a \textsl{hybrid} nodal distribution of nodes, namely we adopt the configuration that combines the uniform distribution of nodes with the Waldron nodes \cite{waldron}. In particular, we consider a $7$-degree uniform distribution of nodes and, in each of the $49$ triangles of the corresponding net, we construct a $10$-degree distribution of Waldron nodes by using the algorithm described in \cite{waldron}. In Figure \ref{fig:bezierwaldronconfiguration}, we show the node distribution, the sparsity pattern of the collocation matrix, and the error on the surface of the approximated solution by using the MS-FD method. 
	The condition number of the collocation matrix for the MS-FD method is about $3.9$ at the first time step, and $2.8$ starting from the second time step onward. In Table \ref{errore:bezier+waldron} we report the mean and maximum absolute error for the various time steps by using the MS-FD approximation. Compared to previous node distributions, the MS-FD methods still show an error of the same magnitude. 
	\begin{table}
		\centering
		\footnotesize
		\begin{tabular}{c c c c c }
			\hline
			Time Step &$E_{mean}^{RBF-FD}$& $E_{max}^{RBF-FD}$  &$E_{mean}^{MBDF2}$&$E_{max}^{MBDF2}$\\
			\hline
			1&2.1679e-3&4.4602e-2&  2.4764e-4&  2.1957e-2\\
			2&2.1658e-3&4.5045e-2&  2.4916e-4&  2.1036e-2\\
			3&2.1554e-3&4.5312e-2&  2.5046e-4&  2.0249e-2\\
			4&2.1393e-3&4.5404e-2&  2.5236e-4&  1.9667e-2\\
			5&2.1182e-3&4.5335e-2&  2.5453e-4&  2.0193e-2\\
			6&2.0927e-3&4.5120e-2&  2.5676e-4&  2.0697e-2\\
			7&2.0638e-3&4.4778e-2&  2.5906e-4&  2.1181e-2\\
			8&2.0314e-3&4.4322e-2&  2.6148e-4&  2.1649e-2\\
			9&1.9965e-3&4.3767e-2&  2.6394e-4&  2.2102e-2\\
			10&1.9598e-3&4.3127e-2&  2.6646e-4&  2.2542e-2\\
			11&1.9220e-3&4.2414e-2&  2.6907e-4&  2.2972e-2\\
			12&1.8839e-3&4.1639e-2&  2.7159e-4&  2.3390e-2\\
			13&1.8455e-3&4.0811e-2&  2.7441e-4&  2.3799e-2\\
			14&1.8071e-3&3.9939e-2&  2.7752e-4&  2.4198e-2\\
			15&1.7695e-3&3.9033e-2&  2.8089e-4&  2.4588e-2\\
			16&1.7322e-3&3.8098e-2&  2.8461e-4&  2.4970e-2\\
			17&1.6956e-3&3.7143e-2&  2.8847e-4&  2.5342e-2\\
			18&1.6598e-3&3.6172e-2&  2.9251e-4&  2.5706e-2\\
			19&1.6246e-3&3.5191e-2&  2.9677e-4&  2.6061e-2\\
			20&1.5908e-3&3.4204e-2&  3.0133e-4&  2.6407e-2\\
			\hline
		\end{tabular}
		\caption{Mean and maximum absolute error for the uniform distribution of nodes combined with Waldron nodes in Figure \ref{fig:bezierwaldronconfiguration} on the various time steps by using the MS-FD approximations.}
		\label{errore:bezier+waldron}
	\end{table}
	
	\subsection{Uniform distribution of nodes combined with Waldron nodes with the addition of nodes lying on special straight lines}
To better capture the behavior of the solution near critical regions, we propose a new strategy that allows us to balance accuracy and stability while increasing resolution in areas of interest.
Starting from the node distribution of Section \ref{sec:bezier+waldron}, to better capture the behavior of the solution near the origin, we intensify it at the boundary, by adding $52$ of equispaced nodes along the lines $x=0$ and $y=0$, from $x=y=0$ up to $x=y=2.5$. Further, to reduce the error near the line $x+y=2$, where the solution shows higher localized variations of the gradient, we add $52$ equally spaced nodes along the lines $x+y = a$, with $a \in \left\{ 1.5,1.8,2,2.1,2.3,2.5,2.7\right\}$. Finally, we add $52$ equally spaced nodes along the lines $x+y = a$, with $a \in \left\{7.1,7.5\right\}$ and we introduce a gap of $0.4$ between the interior points and the boundary $x+y = 8$. This configuration of nodes was inspired by the studies of U. Pettersson et al. \cite{PETTERSSON200882} and S. Milovanović et al. \cite{MILOVANOVIC2020205}, where adaptive node sets with increasing density near the singular line $x+y=2$ and decreasing density toward the line $x+y=8$ were employed.

\begin{figure}
		\centering
		\parbox{0.32\linewidth}{\centering
			\includegraphics[width=\linewidth]{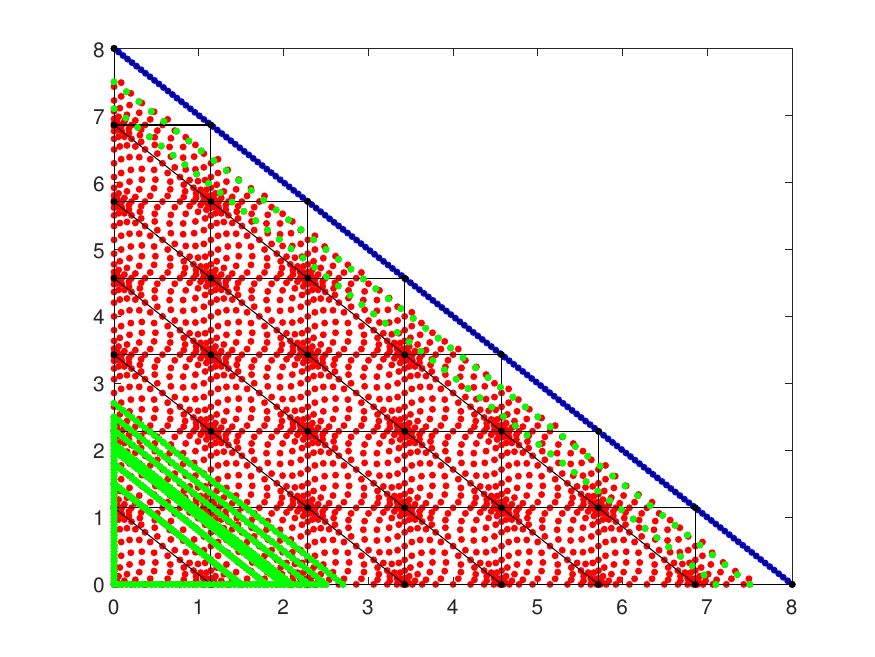}\\[-0.5ex]
			(a)
		}
		\parbox{0.32\linewidth}{\centering
			\includegraphics[width=1\linewidth]{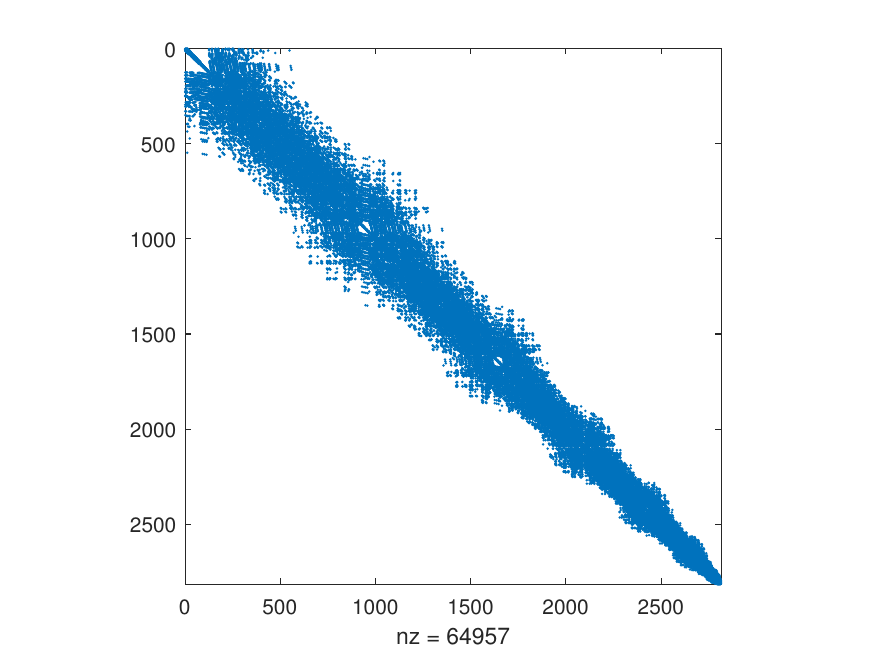}\\[-0.5ex]
			(b)
		}
		\parbox{0.32\linewidth}{\centering
			\includegraphics[width=1\linewidth]{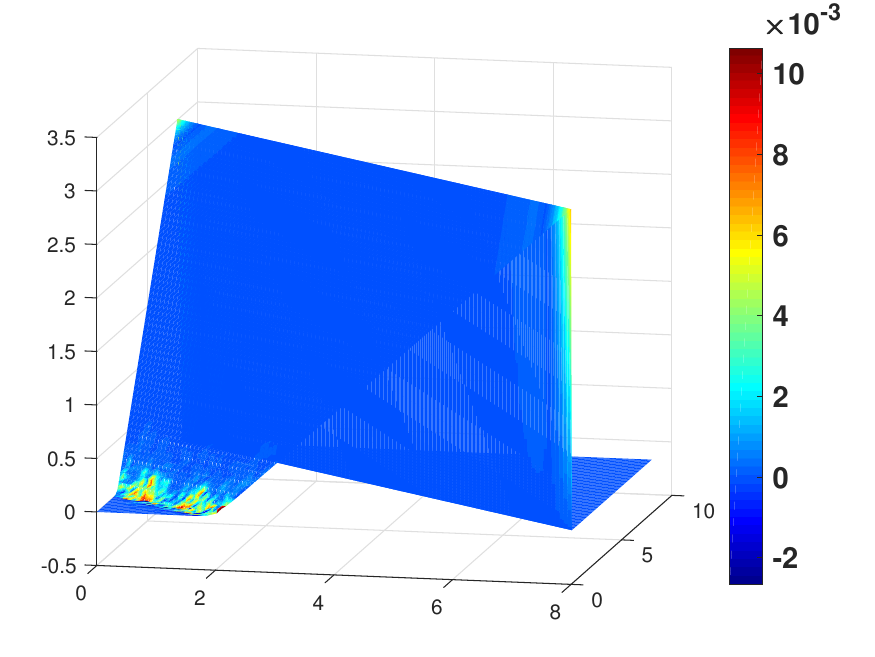}\\[-0.5ex]
			(c)
		}
		
		\caption{(a) Set of $2814$ uniform distribution of nodes combined with Waldron nodes  and straight lines coupled with $153$ boundary points. (b) Sparsity pattern of the collocation matrix. (c) Error on the surface of the approximated solution at the last time step with the MS-FD algorithm.}
		\label{fig:bezierwaldronretteconfiguration}
	\end{figure}

Figure \ref{fig:bezierwaldronretteconfiguration} illustrates the node distribution, the sparsity pattern of the collocation matrix, and the error surface of the approximated solution obtained using the MS-FD method. The condition number of the collocation matrix is approximately $33.3$ at the first time step and $17.9$ from the second time step onward. Table \ref{errore:bezier+waldron+rette} reports the mean and maximum absolute errors at various time steps.

The numerical results show that the mean absolute errors of the MS-FD scheme are on the order of $\approx 10^{-4}$, outperforming the RBF-FD mean absolute errors which are on the order of $\approx 10^{-3}$. As shown in Table \ref{errore:bezier+waldron+rette}, this final nodal distribution is the most accurate among all configurations considered in the paper for the MS-FD scheme, as it reduces the maximum absolute error at certain time steps, while maintaining the mean absolute error within the same order of magnitude.

	\begin{table}
		\centering
		\footnotesize
		\begin{tabular}{c c c c c}
			\hline
			Time Step &$E_{mean}^{RBF-FD}$& $E_{max}^{RBF-FD}$  &$E_{mean}^{MS-FD}$&$E_{max}^{MS-FD}$\\
			\hline
			1  &2.1679e-3&4.4602e-2&  1.7117e-4&  1.2420e-2\\
			2   &2.1658e-3&4.5045e-2&  1.7358e-4&  1.1775e-2\\
			3   &2.1554e-3&4.5312e-2&  1.7485e-4&  1.1208e-2\\
			4   &2.1393e-3&4.5404e-2&  1.7603e-4&  1.0712e-2\\
			5  &2.1182e-3&4.5335e-2&  1.7685e-4&  1.0280e-2\\
			6   &2.0927e-3&4.5120e-2&  1.7738e-4&  9.9035e-3\\
			7   &2.0638e-3&4.4778e-2&  1.7770e-4&  9.5753e-3\\
			8  &2.0314e-3&4.4322e-2&  1.7776e-4&  9.2890e-3\\
			9   &1.9965e-3&4.3767e-2&  1.7770e-4&  9.0387e-3\\
			10 &1.9598e-3&4.3127e-2&  1.7738e-4&  8.8196e-3\\
			11  &1.9220e-3&4.2414e-2&  1.7708e-4&  8.6275e-3\\
			12 &1.8839e-3&4.1639e-2&  1.7676e-4&  8.8742e-3\\
			13 &1.8455e-3&4.0811e-2&  1.7668e-4&  9.2111e-3\\
			14&1.8071e-3&3.9939e-2&  1.7676e-4&  9.5332e-3\\
			15&1.7695e-3&3.9033e-2&  1.7703e-4&  9.8437e-3\\
			16&  1.7322e-3&3.8098e-2&  1.7768e-4&  1.0145e-2\\
			17 &  1.6956e-3&3.7143e-2&  1.7851e-4&  1.0438e-2\\
			18  & 1.6598e-3&3.6172e-2&  1.7956e-4&  1.0724e-2\\
			19 &1.6246e-3&3.5191e-2&  1.8096e-4&  1.1005e-2\\
			20  & 1.5908e-3&3.4204e-2&  1.8274e-4&  1.1280e-2\\
			\hline
		\end{tabular}
		\caption{Mean and maximum absolute error for the uniform distribution of nodes combined with Waldron nodes  and straight lines in Figure \ref{fig:bezierwaldronretteconfiguration} on the various time steps by using the MS-FD approximations.}
		\label{errore:bezier+waldron+rette}
	\end{table}

    \begin{remark}
   In the numerical experiments, we also considered, the Backward Differentiation Formula of third order (BDF3) for time discretization. However, the results did not show any significant improvement compared to the use of the BDF2 scheme. The latter provides stable solutions with a lower computational cost, being a two-step method. Based on the experiments conducted, the use of the BDF3 scheme does not yield substantial benefits in terms of accuracy, as the approximation error remains of the same order of magnitude. The same considerations apply when using local interpolating polynomials of degree higher than $2$. For this reason, we chose to work with degree $2$ polynomials, which are computationally less expensive.
    \end{remark}

	\section{Conclusion}\label{sec7}
    	The main objective of this paper is the introduction of a collocation method based on the multinode Shepard interpolant for the numerical solution of the two-dimensional Black-Scholes equation, with particular focus on the second-order Backward Differentiation Formula (BDF2) for time discretization. The choice of a suitable nodal distribution, a key factor in ensuring solution accuracy, is also addressed throughout the paper.
	We can conclude that the multinode Shepard method, applied via collocation, proves to be effective in approximating the solution of the two-dimensional Black-Scholes equation when combined with implicitly stable time discretisation techniques such as BDF2.
	
	Future work may explore alternative point distributions or the enrichment of the functional space spanned by the multinode cardinal Shepard basis with additional functions to improve accuracy.

\section*{Acknowledgements}

This research has been achieved as part of RITA “Research ITalian network on Approximation”, as part of the UMI group “Teoria dell’Approssimazione e Applicazioni” and was supported by INDAM-GNCS project 2025. The authors are members of the INdAM Research group GNCS.
    
	\bibliographystyle{elsarticle-num}
	\bibliography{bibliografia}
	
\end{document}